\documentstyle{article}
\def\R{{I\!\! R}}
\def\N{{I\!\! N}}
\def\Z{{Z\!\!\! Z}}
\def\kasten{$~~\mbox{\hfil\vrule height6pt width5pt depth-1pt}$ }
\def\C{~\hbox{\vrule width 0.6pt height 6pt depth 0pt \hskip -3.5pt}C}
\newtheorem{Theorem}{Theorem}[section]
\newtheorem{Problem}[Theorem]{Problem}
\newtheorem{Definition}[Theorem]{Definition}
\newtheorem{Proposition}[Theorem]{Proposition}
\newtheorem{Lemma}[Theorem]{Lemma}
\newtheorem{Remark}[Theorem]{Remark}

\begin{document}
\par\noindent
{\large \bf Partly Divisible Probability Distributions}

\bigskip
\par\noindent
S. Albeverio${}^{1,2,3}$, H. Gottschalk${}^1$, and J.-L. Wu${}^{1,2,4}$
\medskip
\par\noindent
${}^1$ Fakult\"at und Institut f\"ur Mathematik der 
Ruhr-Universit\"at Bochum, D-44780 Bochum, Germany
\par\noindent
${}^2$ SFB237 Essen-Bochum-D\"usseldorf, Germany
\par\noindent
${}^3$ BiBoS Research Centre, Bielefeld-Bochum, Germany; and CERFIM,
Locarno, Switzerland
\par\noindent
${}^4$Probability Laboratory, Institute of Applied Mathematics, 
Academia Sinica, Beijing 100080, P R China

\bigskip\par
\begin{abstract}{Given a probability distribution $\mu$ a set 
$\Lambda (\mu)$ of positive real numbers
is introduced, so that $\Lambda (\mu)$ measures the
"divisibility" of $\mu$. The basic properties of $\Lambda
(\mu)$ are described and examples of probability distributions are given, 
which exhibit the existence of a continuum of situations interpolating the
extreme cases of infinitely and minimally divisible
probability distributions.}
\end{abstract}

\bigskip\par\noindent
{\bf Mathematics Subject Classification (1991)}: Primary 42A85; 
\par\noindent
Secondary 60B15, 60E10

\bigskip\par\noindent
{\bf Key Words and Phrases}: admissible probability distributions, convolution semigroups,
divisibility, Fourier transforms.

\section*{Introduction}
In probability theory it is often interesting to know whether a random
variable $X$ can be decomposed into a sum of $n$ independent identical
distributed random variables 
$$X=X_1+\cdots+X_n ~.$$

In terms of the probability distribution $\mu $ of $X$ the above
equation means that
\begin{equation}
\label{0.1eqa}
 \mu = \underbrace{\mu_{1\over n}\ast \cdots \ast \mu_{1\over n}}_{n
~ \mbox{\tiny times}} ,
\end{equation}
where $\mu_{1\over n}$ is the probability distribution of $X _1$.

The notion of infinitely divisible probability distribution
characterizes a special case of such phenomena, where the process of
dividing a probability distribution $\mu $ into $n$-fold convolutions of
probability distributions $\mu_{1\over n}$ can be perpetuated {\it ad
infinitum}.

On the other hand, some random events can only be decomposed into some
finite number of identical elementary random
events, as e.g. the $n$-fold coin-toss which can be decomposed at
maximum into $n$ elementary events, since the single coin-toss is
clearly indecomposable \cite{Lu} (For a recent presentation of "~arithmetic
properties" of probability distributions, see \cite{RUSZ} ).

The aim of this paper is to give a definition of a set 
$\Lambda (\mu )$ of positive real numbers 
which in some sense measures the "divisibility" of a
probability distribution $\mu$. From a dynamical viewpoint, the set
$\Lambda (\mu)$ appears as the index set of some maximal convolution semigroup
of probability measures $ (\mu _t )_{t\in\Lambda(\mu )}$ associated to
the given probability distribution $\mu $. This  generalizes the well-known
notion of convolution semigroup of probability measures $(\mu_t)_{t>0}$
associated with an infinitely divisible measure $\mu$ (see e.g. \cite{BeFo,Hey}).

Let us remark that the idea of considering the set $\Lambda(\mu)$ has 
already been introduced in the classic work \cite{Du}.  In this paper, we 
develop however different techniques for the study of $\Lambda(\mu)$ and
for the construction of illustrative examples for the structure of 
$\Lambda(\mu)$.

\section{Notations and Definitions}
By a \underline{characteristic function} we mean a function $f:\R \to
\C$ with the following properties (c.f. \cite{Lu})
\begin{itemize}
\item[1.] $f$ is continuous;
\item[2.] $f$ is positive-definite;
\item[3.] $f(0)=1$.
\end{itemize}
Let $\mu$ be a probability distribution, i.e. a probability measure on
the real line. Then, the \underline{Fourier transform} $\hat \mu :\R\to\C$ of
$\mu$ is given by
$$\hat \mu (y) := \int _{\R} e^{iyx} \mu(dx) ~ , ~~ y\in \R .$$
Sometimes, we will also use the symbol $\cal F$ for the Fourier
transform and ${\cal F}^{-1}$ for its inverse.
By Bochner's theorem \cite{B} the Fourier transform provides a one to one
correspondence between probability distributions and characteristic
functions. We may thus say a characteristic function
$\hat \mu$ is associated to a probability measure $\mu$ (and vice versa).

We call a probability distribution \underline{admissible}, if the
associated characteristic function maps $\R$ into $\C^*:=\C - \{ 0 \}$.
By ${\cal M}_1^a(\R)$ we denote the collection of all admissible
probability distributions. The sense of this notion is clarified by the
following lemma:

\begin{Lemma}
\label{1.1lem}
Let $\mu \in {\cal M}_1^a(\R)$. Then there exists a uniquely determined
continuous function $\psi :\R\to \C$ with $\psi(0)=0$ and $\hat \mu = e^\psi$.
\end{Lemma}
\noindent {\bf Proof}~Existence: Let $I\subset [0,\infty)$ be the maximal 
subinterval of $[0,\infty)$, for which a function $\psi$ fulfilling the
requirements of the lemma can be defined. Then $I\not = \emptyset$,
since in a right ~neighborhood of zero $\log \hat \mu$ fulfills the
requirements, where $\log$ stands for the principal branch of the complex
logarithm. Suppose $I=[0,a)$ for some $a>0$. Let $y_0\in [0,a)$ such
that for $y\in [y_0,a]$ $\hat \mu (y)$ does not hit the ~semi axis in
the complex plane $\C$ that starts at $0$ and goes through $-\hat \mu
(y_0)$. Thus 
$$\psi (y) = \psi (y_0)+(\psi(y)-\psi(y_0)) =
\psi(y_0)+\log(\hat\mu(y)\hat\mu(y_0)^{-1}).$$
Since the limit of the right hand side for $y\uparrow a$ exists by our
assumptions on $y_0$, $\psi $ can be continuously extended to $[0,a]$,
which is in contradiction with the maximality of $I$.

Similarly we get a contradiction, if we assume $I=[0,a]$ for some $a>0$:
We choose $\epsilon$ small enough, such that $\hat \mu (y)$ for $y\in 
[a,a+\epsilon) $ does not hit the
~semi axis from $0$ through $-\hat \mu (a)$ and we extend $\psi$ continuously
by setting $\psi(y):=\psi(a)+\log(\mu(y)\mu(a)^{-1})$ for $y\in[a,a+\epsilon)$.

Thus we conclude that $I=[0,\infty)$ and a similar consideration
yields that $\psi$ can be well-defined on the closed negative ~semi axis
$I'=(-\infty,0]$, too.

Uniqueness: Suppose $\psi_1$ and $\psi_2$ fulfill the requirements of
the lemma. For $y\in \R$ we have $e^{\psi_1(y)}=e^{\psi_2(y)}$, hence we get
$$\psi_1(y)-\psi_2(y)=2\pi i k(y)~ , ~~\mbox{with}~~k(y)\in \Z ~~
\mbox{for}~~y\in \R.$$
We remark that the left hand side is continuous in $y$, thus the right 
hand side must also  be
continuous, which is possible only if it is constant. Evaluating the
above equation for $y=0$ then yields that the right hand side is equal
to zero for all $y\in\R$. \kasten

For instance all infinitely divisible probability distributions are
admissible \cite{Lu}. The function $\psi$ defined in Lemma \ref{1.1lem}
in the case of an infinitely divisible probability distribution
coincides with the \underline{second characteristic} \cite{Lu} of the
infinitely divisible probability distribution. From now on, we will use
this term also in the general case of admissible probability
distributions.

Clearly, if a characteristic function $\hat \mu$ maps $\R$ into $\C -(-\infty,0]$, 
then the second characteristic of the associated probability distribution $\mu $
equals to $\log \hat \mu$, where $\log $ denotes the main branch of the complex
logarithm.

We are now ready to give the central definition of this article.

\begin{Definition}
\label{1.1def}
For $\mu \in {\cal M}_1^a(\R)$, we define the subset $\Lambda(\mu) $ of positive real
numbers  as follows:
$$\Lambda(\mu) := \{ t>0 : e^{t\psi} ~ \mbox{ is a
characteristic function} \} ~. $$
Here $\psi $ denotes the second characteristic of $\mu$.
\end{Definition}

We remark that \cite{Du}(see page 45) contains a corresponding definition 
of the set $\Lambda(\mu)$, however without giving precise conditions on
the class of probability distributions $\mu$ considered. 

In our terminology, the statement of Schoenberg's theorem \cite{BeFo} can
be formulated as follows: $\mu $ is infinitely divisible if and only if
$\mu $ is admissible and $\Lambda (\mu ) = (0,\infty)$. In this paper we
want to ask:

\begin{Problem} 
\label{1.1prob}
What kind of set is $\Lambda (\mu)$ if $\mu $ is not
infinitely divisible?
\end{Problem}

The rest of the paper is organized as follows.
In Section 2, we collect some general properties of the set
$\Lambda(\mu)$ for $\mu \in {\cal M}_1^a(\R)$. 

In Section 3, we explain what information concerning 
the divisibility of $\mu$ we can extract from $\Lambda(\mu)$.

Finally in Section 4,  we present concrete examples for $\mu$ including the case
$\Lambda (\mu)=\N$ and the case where $\Lambda(\mu)$ includes a semi-bounded
interval but not $(0,\infty )$. Since in the latter case the lower 
bound of the semi-bounded
interval can go to zero, our set $\Lambda(\mu)$ gives a "calibration" of
divisibility of  admissible probability distributions $\mu$, 
interpolating between the extreme cases $\Lambda(\mu) = (0,\infty)$
and $\Lambda(\mu) = \N$.

\section{Some Properties of $\Lambda(\mu)$}

In order to derive some basic properties for the set $\Lambda(\mu)$, let
us first recall the following identity for the Fourier transform of the
convolution of two probability measures $\mu$ and $\nu$:
\begin{equation}
\label{2.1eqa}
{\cal F} (\mu \ast \nu ) = {\cal F} (\mu) ~ {\cal F} (\nu)
\end{equation}
In particular, the product of two characteristic functions is also a
characteristic function. This gives rise to the following

\begin{Proposition}
\label{2.1prop}
Let $\mu \in {\cal M}_1^a(\R)$. Then $\Lambda (\mu)$ is a semigroup
containing $1$.
Furthermore, $\Lambda (\mu)$ is closed in $(0,\infty)$.
\end{Proposition}
\noindent {\bf Proof} Let $\psi$ denote the second characteristic of
$\mu$. Since $\hat{\mu}=e^{1\cdot\psi}, \, 1\in\Lambda(\mu)$.

Suppose that $s,t \in \Lambda(\mu)$, then $e^{s\psi}$ and
$e^{t\psi}$ are characteristic functions. It follows that 
$e^{s\psi}e^{t\psi}=e^{(s+t)\psi}$ is also a characteristic function and
thus $s+t\in \Lambda(\mu)$. The conclusion that $1\in\Lambda (\mu)$ is 
by the definition of $\Lambda (\mu)$.

Furthermore, we choose a sequence $\{ t_n\}_{n\in \N}\subset
\Lambda(\mu)$ and we assume that $t_n \to t$, where $t\in (0,\infty)$.
Clearly $e^{t\psi}$ is continuous and fulfills $e^{t\psi (0)} = 1$.
Since the pointwise limit of positive definite functions is a positive
definite function and $e^{t_n\psi}\to e^{t\psi}$ pointwise, $e^{t\psi}$
is also positive definite and thus has all properties of a characteristic 
function. Therefore $t\in\Lambda(\mu)$. \kasten

Since $1\in\Lambda(\mu)$ for any $\mu \in {\cal M}_1^a(\R)$, an immediate
consequence of Proposition \ref{2.1prop} is that $\N \subset \Lambda (\mu)$.
From now on we call an admissible probability measure $\mu $ on $\R$ with 
$\Lambda(\mu) = \N$ \underline{minimally divisible}, since $\N$ is the
smallest (closed) semigoup in $(0,\infty)$ including $1$. This situation
clearly is opposite to the infinitely divisible situation, where
$\Lambda(\mu) =(0,\infty)$, i.e., $\Lambda (\mu)$ is maximal (namely, infinite
divisible probability measures are "maximally divisible").

For $\mu \in {\cal M}_1^a(\R)$ with second characteristic $\psi $ 
and $t\in \Lambda(\mu)$ there exists a
uniquely determined measure $\mu_t$ with characteristic function $\hat
\mu _t := e^{t\psi}$. Applying the inverse Fourier transform ${\cal F}^{-1}$ 
to equation (\ref{2.1eqa}), we get that 
the family of probability distributions $(\mu_t)_{t\in\Lambda(\mu)}$ forms a
convolution semigroup of probability measures, i.e. for $s,t \in
\Lambda(\mu)$, we obtain $\mu_s\ast\mu_t=\mu_{s+t}$.

Thus, $\Lambda(\mu)$ can be viewed as an index-set for the "time-index" $t$
of the convolution semigroup of probability measures $(\mu_t)_{t\in
\Lambda(\mu)}$.
From the definition of $\Lambda (\mu)$, it is clear that $(\mu_t)_{t\in
\Lambda(\mu)}$ is the maximal convolution semigroup which can be
associated with the admissible probability distribution $\mu$ by the
requirement that $\mu = \mu_1$.

Furthermore, for $\mu \in {\cal M}_1^a(\R)$ with second characteristic $\psi$ 
and $t\in \Lambda(\mu)$, we
also have $\mu _t \in {\cal M}_1^a(\R)$ and the second characteristic of
$\mu_t$ is clearly $t\psi$. From Definition \ref{1.1def} we can then get the
following relation
\begin{equation}
\label{2.2eqa}
\Lambda(\mu_t) = {s\over t} \Lambda (\mu_s) ~~~ s,t\in \Lambda(\mu).
\end{equation}
Equation (\ref{2.2eqa}) gives a simple scaling property for the family
$\{\Lambda(\mu_t)\}_{t\in \Lambda(\mu)}$.

The following proposition  investigates some properties of general
closed semigroups $S\subset (0,\infty)$ including $1$. Here $S$
does not necessarily have to be associated with some admissible probability
distribution.

\begin{Proposition}
\label{2.2prop}
Let $S\subset (0,\infty)$ be a closed semigroup with $1\in S$.
\begin{itemize}
\item[1.] Either $S=(0,\infty)$ or $S\subset[\lambda,\infty)$ for some
$\lambda >0$.
\item[2.] If the interior of $S$ is non-empty, then $S\supset
[\lambda,\infty)$ for some $\lambda, ~0<\lambda <~\infty$.
\end{itemize}
\end{Proposition}
\noindent {\bf Proof} {\it 1.} Suppose that $\not \exists \lambda >0 $
such that $S\subset [\lambda,\infty)$. Then there is a sequence
$\{t_n\}_{n\in\N}\subset S$ with $t_n \to 0$. Now we fix an arbitrary
$t\in (0,\infty)$ and we choose a subsequence $\{s_n\}_{n\in\N}$, say, of
$\{t_n\}_{n\in\N}$ such that $s_n<t$ and $s_n<{1\over n}$ for all $n\in
\N$. Let $m_n$ be the largest integer with $m_ns_n\leq t$. By the
semigroup property of $S$ we have $m_ns_n\in S$. Since by construction,
$|m_ns_n-t|<{1\over n}$,  we then get $m_ns_n\to t$ and the closedness of $S$
now implies $t\in S$. Thus $S=(0,\infty)$.

{\it 2.} If the open interior of $S$ is non-empty, there exists an 
non-empty open ~interval $(a,b)\subset S$. Let $n\in \N$ be sufficiently large
such that $n(b-a) > 1$. By the semigroup property, $(na,nb)\subset S$.
Thus, for $t>na$, there is $t'\in (na,nb)$ and a number $m\in\N\cup\{0\}$
such that $t=t'+m$. Since $\N\subset S$ and $t'\in S$, by the
semigroup property, we get $t=t'+m\in S$. Since $S$ is closed, we 
thus proved $S\supset [na,\infty)$. \kasten

\section{The Connection with Divisibility}

Motivated by the result of Proposition \ref{2.2prop}, we define for $\mu \in 
{\cal M}_1^a(\R) $ the following two parameters, which give a
rough characterization of the set $\Lambda(\mu)$:
\begin{eqnarray}
\label{2.3eqa}
\lambda_0(\mu) &:=& \inf \Lambda(\mu) \nonumber \\
\lambda_1(\mu) &:=& \left\{ \begin{array}{l}\infty ~,~~\mbox{ if the
open interior of $\Lambda(\mu)$ is empty }\\ \inf \{ \lambda >0 :
\Lambda(\mu)\supset [\lambda,\infty) \} ~,~~\mbox{otherwise} \end{array}
\right.
\end{eqnarray}
To connect these two parameters with the problem of divisibility,
as it was posed in equation (\ref{0.1eqa}), we introduce the following
notions. First of all, every non-infinitely divisible probability distribution
$\mu \in {\cal M}_1^a$ is called \underline{partly divisible}. A partly divisible 
probability distribution $\mu$ is also called \underline{$n$-divisible} if $n$
is the largest integer such that equation (\ref{0.1eqa}) has a solution.
Furthermore, we say that $\mu\in{\cal M}_1^a(\R)$ is 
\underline{divisible up to $m$}, if ~equation (\ref{0.1eqa}) can be solved
for $n=1,2,\ldots , m$ and $m$ is the maximal integer with this
property. Here we include infinitely divisible probability distributions
$\mu\in{\cal M}^a_1(\R)$ by saying that $\mu$ is $\infty$-divisible or divisible up to
$\infty$ respectively, which is clearly the same (c.f. Prop.
\ref{2.2prop}.{\it 1.}). Finally, we adopt the
conventions that ${1\over \infty}= 0$ and ${1\over 0}=\infty$.

For an admissible probability distribution $\mu$ which is $n$-divisible
and divisible up to $m$, we set $\tilde \lambda_0(\mu):={1\over n}$ and
$\tilde \lambda_1(\mu):={1\over m}$. Using the fact that any solution
$\mu_{1\over q}$ of an equation of the type (\ref{0.1eqa})
by equation (\ref{2.1eqa}) must be admissible and also using  the scaling
property (\ref{2.2eqa}), we conclude that
$\tilde \lambda_0(\mu) \in \Lambda(\mu)$ and $1, {1\over 2}, \ldots , 
{1\over m}= \tilde \lambda_1(\mu) \in \Lambda(\mu)$. Consequently, we
have the following inequalities:
$$0\leq \lambda_0(\mu)\leq\tilde \lambda_0(\mu)\leq \tilde
\lambda_1(\mu) \leq [\lambda_1(\mu)^{-1}]^{-1}\leq \infty. $$
Here $[r]$ stands for the integer part of $r\in [0,\infty]$ with
the convention that $[\infty]=\infty$.
Thus the parameters $\lambda_0(\mu)^{-1}$ and $[\lambda_1(\mu)^{-1}]$
give, respectively, the
upper and lower bounds for the $n$-divisibility and divisibility up to
$m$ of a probability distribution $\mu \in {\cal M}_1^a(\R)$. To explain
further the connection between the parameters $\lambda_j(\mu)$ and 
$\tilde \lambda_j(\mu), ~ j=0,1$, let us introduce the notation
$\mu^{*n}$ for the $n$-fold convolution of a probability measure with
itself. With this notation, equation (\ref{0.1eqa}) reads $\mu =
(\mu_{1\over n})^{*n}$. 

In order to generalize this equation to non-
integer $n$, for $\mu\in {\cal M}_1^a(\R)$ with second characteristic $\psi$
and $t \in \Lambda(\mu)$ we define $\mu^{*t}:=\mu_t={\cal F}^{-
1}(e^{t\psi})$. By equation (\ref{2.1eqa}) we see that this coincides
with the definition given before in the case that $t\in\N$. We now
generalize the problem posed by equation (\ref{0.1eqa}) to the problem
whether for a given $\beta \in (0,\infty)$ there exists an admissible
measure $\mu_{1\over \beta}$ which solves
\begin{equation}
\label{2.4eqa}
\mu = (\mu_{1\over \beta})^{*\beta} ~.
\end{equation}
By the scaling property (\ref{2.2eqa}), a solution of equation
(\ref{2.4eqa}) exists if and only if ${1\over \beta} \in \Lambda(\mu)$.
If we call a solution of equation (\ref{0.1eqa}) an \underline{integer
root} of $\mu$, then by the above motivations it seems natural to speak of a 
solution $\mu_{1\over \beta}$ for general $\beta \in (0,\infty)$ of
(\ref{2.4eqa}) as a \underline{real root} of $\mu$.

Clearly for a partly divisible $\mu \in {\cal M}_1^a(\R)$,
$\beta=\lambda_0(\mu)^{-1}$ is the largest number for which there exists a solution 
of equation (\ref{2.4eqa}). We thus call that $\mu \in {\cal M}_1^a(\R)$
\underline{$\lambda_0(\mu)^{-1}$-real divisible}. ~Analogously $\mu$
is called \underline{real divisible up to $\lambda_1(\mu)^{-1}$},
since $\beta':=\lambda_1(\mu)^{-1}$ is the largest number, such that equation
(\ref{2.4eqa}) can be solved for all $\beta \in (0,\beta')$. Here again
the cases of $\infty$-real divisibility and real divisibility up to $\infty$
coincide in the case of infinite divisibility. 

As we have seen above, the  $\lambda_j(\mu),~j=0,1$ 
can be interpreted as parameters which measure 
the divisibility of a probability distribution $\mu \in {\cal
M}_1^a(\R)$. Applying the scaling property (\ref{2.2eqa}) to the definition of
$\lambda_j(\mu)$ we immediately get
$$ \lambda_j(\mu_t) = {1\over t}\lambda_j(\mu) ~~\mbox{for}~j=0,1,~
t\in \Lambda(\mu).$$
Since $\N \subset \Lambda(\mu)$, we can deduce that for $\mu \in {\cal M}_1^a(\R)$ 
$\lambda_0(\mu)$ can be 
arbitrarily close to zero. In the
following section we will construct an admissible probability distribution
$\mu$ such that $\Lambda(\mu)$ has non-empty open interior. Therefore, there
exist examples for admissible probability distributions $\mu$ with
$\lambda_1(\mu)$ arbitrarily close to zero.

Thus, from the both viewpoints -- the real divisibility one and the one of
real divisibility up to some positive real number -- there exist examples of
partly divisible ~probability measures which are arbitrarily close to infinite
divisible probability measures. This is what we called in the introduction a
"calibration" of divisibility interpolating between the two extreme cases
$\Lambda(\mu)=(0,\infty)$ and $\Lambda(\mu)=\N$.

\begin{Remark}
\label{3.1lem}
We would like to point out that it is natural to extend our
considerations to the case of probability distributions on locally
compact abelian groups \cite{BeFo,Hey}. However, since Lemma \ref{1.1lem} depends
crucially on the topological structure of the real line $\R$, the
notion of admissibility in the framework of locally compact abelian 
groups has to be revised. We will investigate this problem in the
forthcoming work \cite{AGW}. 
\end{Remark}

\section{Examples}
In this section we present two examples: one example showing the existence of
minimally divisible probability distributions and the other one presenting 
a partly divisible probability distribution $\mu$ with an associated set 
$\Lambda(\mu)$ which has non-empty open interior and therefore
$0<\lambda_1(\mu)<\infty$.

This is certainly no more than an initial step towards the mathematical
problem of determining all closed semigroups $S\subset (0,\infty)$
including 1 which can  ~~occur as index sets of convolution semigroups of
probability measures $(\mu_t)_{t\in\Lambda(\mu)}$ associated with some
$\mu\in {\cal M}_1^a(\R)$ such that $\mu_1=\mu$.

Let us remark that \cite{Du}(see page 45) provides two examples where 
$\Lambda(\mu)=\{1,\frac{3}{2},2,\frac{5}{2},\cdots\}$ for
a special (non-admissible) finite positive measure $\mu$ and respectively 
where $\Lambda(\mu)=[\frac{1}{2}, \infty)$ for a specific (admissible) finite
positive measure $\mu$.

\subsection{A Minimally Divisible Probability Distribution}
Let $\theta$ be a given probability distribution and let $\alpha \in (0,1)$.
We set 
$$\nu := (1+\alpha)^{-1} (\delta_0+\alpha \theta)~.$$
 Then $\nu$ is
admissible and $\hat \nu$ maps $\R$ into $\C - (-
\infty,0]$, since $\hat \nu = (1+\alpha)^{-1}(1+\alpha \hat \theta)$ and
$|\hat\theta |\leq 1$. Therefore the second characteristic of $\nu$ is simply
$\log \hat \nu$. For $t\in (0,\infty)$, we get that $e^{t\log\nu} =
\nu^t$. Using the series expansion formula for the main branch of the
complex function $z\in \C \mapsto (1+z)^t\in\C$ at $0$ for $|z|<1$,
we calculate further and get
\begin{equation}
\label{4.1eqa}
\hat \nu ^t = {1\over (1+\alpha)^t}\sum_{l=0}^\infty \left(
\begin{array}{c} t \\ l \end{array} \right) \alpha ^l \hat \theta ^l~.
\end{equation}
The binomial coefficient $ \left(\begin{array}{c} t \\ l \end{array} \right)$
for $t\in (0,\infty)$ is defined as follows:
\begin{equation}
\label{3.2eqa}
\left( \begin{array}{c} t \\ l \end{array} \right) := { t (t-1) \cdots (t-
l+1)\over l!} ~, ~~\mbox{for}~l\in \N,
\end{equation}
and $\left( \begin{array}{c} t \\ l \end{array} \right)
:=1$ for $l=0$. We do not know so far, whether for an arbitrarily fixed $t\in
(0,\infty)$, $\hat \nu ^t$ is a characteristic function. But certainly
$\hat \nu ^t$ is a tempered distribution \cite{RS}. Therefore, in
order to calculate the inverse Fourier transforms of $\hat \nu ^t$ we can
apply the theory of the Fourier transform of tempered distributions
which extends the theory of the Fourier transform of probability
measures \cite{RS}. We also denote the extended Fourier transform and
its inverse by ${\cal F}$ and ${\cal F}^{-1}$, respectively.
We recall that equation (\ref{2.1eqa}) holds also for the extended
Fourier transform (if at least one of both sides of (\ref{2.1eqa}) is well-defined)
\cite{RS}. Taking the inverse Fourier transform of the left
hand side of (\ref{4.1eqa}), by the linearity and weak continuity of
${\cal F}^{-1}$ \cite{RS} as well as the uniform
convergence of the right hand side of (\ref{4.1eqa}), we
get that
\begin{eqnarray}
\label{4.3eqa}
{\cal F}^{-1}(\nu^t) &=& {1\over (1+\alpha)^t} \sum_{l=1}^\infty 
\left( \begin{array}{c} t \\ l \end{array} \right) \alpha^l {\cal F}^{-
1}(\hat \theta^l) \nonumber \\
& = & {1\over (1+\alpha)^t} \sum_{l=1}^\infty 
\left( \begin{array}{c} t \\ l \end{array} \right) \alpha^l \theta
^{*l}~,
\end{eqnarray}
where we have inductively applied equation (\ref{2.1eqa}) and the 
convention that $\theta ^{*0}=\delta _0$. Taking now
$\theta = \delta_1$, we obtain
\begin{equation}
\label{4.4eqa}
{\cal F}^{-1}(\hat \nu ^t) = {1\over (1+\alpha)^t} \sum_{l=1}^\infty 
\left( \begin{array}{c} t \\ l \end{array} \right) \alpha^l \delta_l ~,
\end{equation}
since $(\delta_1)^{*l} = \delta _l$ for $l\in \N\cup\{ 0 \}$.

For $t\in \N$, all the binomial coefficients are positive for $l\leq t$ and
zero for $t<l$. Thus, the right hand side of (\ref{4.4eqa}) is a
positive measure with total mass 1 and therefore a probability
distribution. On the other hand, if $t$ is not an integer, the sign of the
binomial coefficients is positive for $l=0,\ldots ,[t]$. For
$l=[t]+1, [t]+2, \ldots ,$ the sign of the binomial coefficients 
alternates from $+$ to $-$ by starting with the $+$ sign. Therefore, the
right hand side of (\ref{4.4eqa}) is a signed measure \cite{H} with non-
zero negative part. Consequently, for $t>0,~t\not \in \N$ we have $t
\not \in \Lambda(\nu)$. Thus $\Lambda(\nu) = \N$.

\subsection{A Partly Divisible Probability Distribution which is Real
Divisible Up To $\rho >0$}

We take again, as in the above example, that $\nu := (1+\alpha)^{-1}
(\delta_0+\alpha \delta_1)$ for $\alpha \in (0,1)$. By $\gamma _t$ we
denote the Cauchy distribution with parameter $t\in (0,\infty)$, i.e.
$$\gamma_t(dx) = {1\over \pi}{t\over x^2+t^2} dx~.$$
It is well-known (see e.g. \cite{Lu,BeFo,Hey}) that $\gamma_t$ is infinitely
divisible and has the characteristic function $\hat \gamma_t (y)=e^{-
t|y|}, ~ y\in \R$. 

We consider the measure $\mu := \gamma_1\ast \nu$, which is clearly 
admissible and has second characteristic $\log \hat \mu=-|y|+\log\hat\nu$. 
Now we want to calculate ${\cal F}^{-1}(e^{t\log\hat\mu}) 
= {\cal F}^{-1}(\hat \gamma_1^t\hat \nu ^t)$ for $t\in (0,\infty)$
in the sense of tempered distributions. Applying equation 
(\ref{2.1eqa}) to the Fourier transform of tempered distributions, we get
that
\begin{eqnarray}
\label{4.5eqa}
{\cal F}^{-1}(\hat \mu^t) &=& \gamma_t \ast {\cal F}^{-1}(\hat \nu^t) 
\nonumber \\
&=& {1\over (1+\alpha)^t}\sum_{l=0}^\infty \left( \begin{array}{c} t 
\\ l \end{array} \right) \alpha ^l\gamma_t\ast\delta_l ~,
\end{eqnarray}
where we have applied formula (\ref{4.4eqa}). Taking into account that
$$\gamma_t\ast\delta_l(dx)={1\over\pi}{t\over (x-l)^2+t^2}dx$$ 
for $l\in\N\cup \{0\}$, we can interpret the right hand side of (\ref{4.5eqa}) as
a possibly signed measure with the density (~with respect to. $dx$)
\begin{equation}
\label{4.6eqa}
f_t(x) := {1\over \pi (1+\alpha)^t}\sum_{l=0}^\infty 
\left( \begin{array}{c} t \\ l \end{array} \right) \alpha^l {t\over (x-
l)^2+t^2}~,~~x\in\R,~ t\in (0,\infty).
\end{equation}

Clearly, for fixed $t\in(0,\infty)$, ${\cal F}^{-1}(\hat \mu ^t)$ is a probability 
distribution if and only if $f_t$ is a non-negative function (the fact
that then the total mass of ${\cal F}^{-1}(\hat \mu ^t)$ is $1$ is a consequence of
$\hat \mu^t(0)=1$). Thus $t\in
\Lambda(\mu) \Leftrightarrow f_t\geq 0$.

Now we choose $t\in [{\sqrt \alpha\over 1-\alpha},\infty)$. We want to
prove that for a such number $t$
\begin{equation}
\label{4.7}
\alpha ^l {t\over (x-l)^2+t^2}\geq \alpha ^{l+1} {t\over (x-
(l+1))^2+t^2} \, , ~~\forall x\in \R ~.
\end{equation}
Substituting $y:=x-(l+1)$, we get that this is equivalent to
$$y^2+{2\alpha\over1-\alpha} y + (t^2-{\alpha \over 1-\alpha}) \geq 0 
\, ,~~\forall y\in \R ~.$$
But the latter inequality holds, since the discriminant of the quadratic 
polynomial on the left hand
side is non-positive by our choice $t\geq {\sqrt \alpha \over 1-
\alpha}$.

Furthermore, for $l=[t]+1,[t]+2,\ldots$ we observe that the ~binomial coefficient
$\left( \begin{array}{c} t \\ l \end{array} \right)$ is either zero (if
$t$ is integer) or the sign of the binomial coefficient alternates
starting with the positive sign. For non-integer $t$ and
$l>[t]$, we get the following estimate
$$ \left| {\left( \begin{array}{c} t \\ l \end{array} \right)\over 
\left( \begin{array}{c} t \\ l+1 \end{array} \right)}\right| = \left|
{l+1\over l-t} \right| >1 ~. $$
We conclude that for $x\in \R $ and $t\in [{\sqrt\alpha \over 1-
\alpha},\infty)$, $f_t(x)\geq 0$ since ({\it 1}) the summands on the right
hand side of equation (\ref{4.6eqa}) are positive for $l=0,\ldots, [t]$,
and ({\it 2}) the remaining summands are either zero (if $t\in \N$) or have 
alternating sign starting with the positive 
sign and the absolute value of these summands decreases monotonically.

We have thus proved that $[{\sqrt \alpha\over 1-\alpha},\infty)\subset
\Lambda(\mu)$.

Finally, it remains to prove that $\mu$ is partly divisible, i.e. non-
infinitely divisible. This is equivalent to $f_t(x)<0$ for some $t\in
(0,\infty),~ x\in \R$. To verify this, we let $x=2$ and we calculate
\begin{eqnarray*}
f_t(2)&=&{1\over \pi (1+\alpha)^t}\left( {\alpha^2\over 2} {t^2(t-1)\over
t^2} + \sum_{l=0,l\not = 2}^\infty \left( \begin{array}{c} t \\ l
 \end{array} \right) \alpha ^l {t\over (2-l)^2+t^2} \right) \\
&<& {1\over \pi (1+\alpha)^t} \left( {\alpha ^2\over 2}(t-
1)+t\sum_{l=0,l\not = 2}^\infty \left| \left( \begin{array}{c} t \\ l
 \end{array} \right)\right| \alpha^l\right) \longrightarrow
 - {\alpha^2\over 2\pi} 
\end{eqnarray*}
 as $ t \to +0$, which implies that $f_t(2)<0$ if $t>0$ is small enough.

{\bf Acknowledgements} The financial support of D.F.G. (SFB 237) is
gratefully acknowledged. The authors also would like to thank the
referee for the careful reading of the paper.

\end{document}